\documentclass[a4paper]{amsart}
\usepackage{graphics}
\usepackage{amsmath}
\usepackage{latexsym}
\usepackage{amssymb}
\usepackage[all]{xy}
\xyoption{matrix}
\xyoption{arrow}

\begin{document}

\renewcommand{\th}{\operatorname{th}\nolimits}
\newcommand{\rej}{\operatorname{rej}\nolimits}
\newcommand{\extto}{\xrightarrow}
\renewcommand{\mod}{\operatorname{mod}\nolimits}
\newcommand{\Sub}{\operatorname{Sub}\nolimits}
\newcommand{\ind}{\operatorname{ind}\nolimits}
\newcommand{\Fac}{\operatorname{Fac}\nolimits}
\newcommand{\add}{\operatorname{add}\nolimits}
\newcommand{\Hom}{\operatorname{Hom}\nolimits}
\newcommand{\Rad}{\operatorname{Rad}\nolimits}
\newcommand{\RHom}{\operatorname{RHom}\nolimits}
\newcommand{\uHom}{\operatorname{\underline{Hom}}\nolimits}
\newcommand{\End}{\operatorname{End}\nolimits}
\renewcommand{\Im}{\operatorname{Im}\nolimits}
\newcommand{\Ker}{\operatorname{Ker}\nolimits}
\newcommand{\Coker}{\operatorname{Coker}\nolimits}
\newcommand{\Ext}{\operatorname{Ext}\nolimits}
\newcommand{\op}{{\operatorname{op}}}
\newcommand{\Ab}{\operatorname{Ab}\nolimits}
\newcommand{\id}{\operatorname{id}\nolimits}
\newcommand{\pd}{\operatorname{pd}\nolimits}
\newcommand{\A}{\operatorname{\mathcal A}\nolimits}
\newcommand{\C}{\operatorname{\mathcal C}\nolimits}
\newcommand{\D}{\operatorname{\mathcal D}\nolimits}
\newcommand{\X}{\operatorname{\mathcal X}\nolimits}
\newcommand{\Y}{\operatorname{\mathcal Y}\nolimits}
\newcommand{\F}{\operatorname{\mathcal F}\nolimits}
\newcommand{\Z}{\operatorname{\mathbb Z}\nolimits}
\renewcommand{\P}{\operatorname{\mathcal P}\nolimits}
\newcommand{\T}{\operatorname{\mathcal T}\nolimits}
\newcommand{\G}{\Gamma}
\renewcommand{\L}{\Lambda}
\renewcommand{\r}{\operatorname{\underline{r}}\nolimits}
\newtheorem{lem}{Lemma}[section]
\newtheorem{prop}[lem]{Proposition}
\newtheorem{cor}[lem]{Corollary}
\newtheorem{thm}[lem]{Theorem}
\newtheorem*{thmA}{Theorem A}
\newtheorem*{thmB}{Theorem B}
\newtheorem{rem}[lem]{Remark}
\newtheorem{defin}[lem]{Definition}

%Egne definisjoner%

\title[Cluster-tilted algebras]{Cluster-tilted algebras}

\author[Buan]{Aslak Bakke Buan}
\address{Institutt for matematiske fag\\
Norges teknisk-naturvitenskapelige universitet\\
N-7491 Trondheim\\
Norway}
\email{aslakb@math.ntnu.no}

\author[Marsh]{Robert J. Marsh}
\address{Department of Mathematics \\
University of Leicester \\
University Road \\
Leicester LE1 7RH \\
England
}
\email{rjm25@mcs.le.ac.uk}

\author[Reiten]{Idun Reiten}
\address{Institutt for matematiske fag\\
Norges teknisk-naturvitenskapelige universitet\\
N-7491 Trondheim\\
Norway}
\email{idunr@math.ntnu.no}


\begin{abstract}

We introduce a new class of algebras, which we call \emph{cluster-tilted}.
They are by definition the endomorphism algebras of
tilting objects in a cluster category. We show that their representation 
theory is very close to the representation theory of hereditary algebras.
As an application of this, we prove a generalised version of so-called
APR-tilting.
\end{abstract}

\thanks{Aslak Bakke Buan was supported by a grant from the Norwegian Research Council.
Robert Marsh was supported by a Leverhulme Fellowship and
by the Norwegian Research Council}

\maketitle

\section*{Introduction}

In~\cite{fz} Fomin and Zelevinsky introduced a class of algebras which 
they called cluster algebras.
There are interesting connections to their theory in many directions (see~\cite{fz2}), amongst them 
to tilting modules over hereditary algebras~\cite{mrz}. Motivated by~\cite{mrz},
tilting theory in a particular factor category of the bounded derived category 
$D^b(\mod H)$ 
of the finitely generated modules over a finite dimensional hereditary algebra $H$ over
a field $k$ was investigated in~\cite{bmrrt}, along with the relationship to tilting 
theory for hereditary 
algebras. In the triangulated category  $D^b(\mod H)$ denote by $[1]$ the suspension functor. Since $H$
has finite global dimension, $D^b(\mod H)$ has AR-triangles~\cite{hap}. Denote by $\tau$
the corresponding translation functor, which is an autoequivalence. Denoting by $F$ the composition
$\tau^{-1}[1]$, the cluster category $\C$ was defined in~\cite{bmrrt} as the factor category  $D^b(\mod H) / F$.

Tilting modules play an important role in the representation theory of finite dimensional algebras, and
the tilted algebras which by definition are the algebras of 
the form $\End_H(T)^{\op}$ for a tilting module $T$
over a hereditary algebra $H$, form a central class of algebras. Their module theory is to a 
large extent determined by the module theory of the hereditary algebra $H$. This motivates 
investigating the cluster-tilted algebras more closely, which are those of the form $\End_{\C}(T)^{\op}$,
where $T$ is a tilting object in $\C$, and the relationship between
$\C$ and the module theory of these algebras. In some sense the relationship is even closer
than is the case for tilted algebras. The following first main result is proved
in Section 2.

\begin{thmA}

If $T$ is a tilting object in $\C$, then
$\Hom_{\C}(T,\ )$ induces an equivalence $\C / \add (\tau T) \to \mod \End_{\C}(T)^{\op}$.
\end{thmA}

Interesting consequences are that for finite representation type
the cluster-tilted algebras all have the same number of indecomposable modules as the 
hereditary algebra $H$ we started with, and that we get 
a nice description of the AR-quivers of cluster-tilted algebras.

Also the notion of almost complete tilting object was investigated in
cluster categories in~\cite{bmrrt}, and it was shown that there are 
always exactly two complements (see Section 1 for definitions). For an almost 
complete tilting object $\overline{T}$, let $M$ and $M^{\ast}$ be the two complements,
and let $T= \overline{T} \oplus M$ and $T' = \overline{T} \oplus M^{\ast}$.
As an application of Theorem~A we prove a close connection between the 
modules over $\G = \End_{\C}(T)^{\op}$ and 
$\G' = \End_{\C}(T')^{\op}$. 
To understand this connection better is of interest for the relationship to 
cluster algebras. If $S_M$ and $S_{M^{\ast}}$ denote the simple
tops of the $\G$-module $\Hom_{\C}(T,M)$ and the $\G'$-module
$\Hom_{\C}(T',M^{\ast})$, respectively, we have the 
following, which answers a conjecture in~\cite{bmrrt}.

\begin{thmB} \label{ThmB}

With the above notation the categories $\mod \G / \add S_M$ and \linebreak[4] 
$\mod \G' / \add S_{M^{\ast}}$
are equivalent.
\end{thmB}

A relevant model for our result from tilting theory is the APR-tilting
modules, which by definition, for a hereditary (basic) algebra $H$
are of the form $T = P \oplus \tau^{-1} S$, where $H = P \oplus S$, and $S$
is simple projective. Then $T$ and $H$ are two completions of the almost complete 
tilting module $P$, and are actually also the two completions of
the almost complete tilting object $P$ in $\C$. Thus, Theorem~B can be regarded
as a generalisation of APR-tilting.

P.~Caldero, F.~Chapoton and R.~Schiffler~\cite{ccs} have recently
associated an algebra to each cluster $C$ in a cluster algebra of simply-laced
Dynkin type, giving a definition via the combinatorics of the
cluster algebra. It is conjectured~\cite[9.2]{bmrrt} that this algebra
coincides with the cluster-tilted algebra associated to the tilting object
in the cluster category $\C$ corresponding to the cluster $C$. They have
given an interesting
geometric description of the module category of this algebra in type
$A_n$, and their approach enables them to generalise the denominator theorem
of Fomin and Zelevinsky~\cite[1.9]{fz2} to an arbitrary cluster in type
$A_n$. Our main aims in this paper are to show that the representation theory
of cluster-tilted algebras is very close to the representation theory of
hereditary algebras, and that they satisfy a generalised version of
APR-tilting.

Theorem~B is proved in Section 4. In the second section we prove Theorem~A,
and in the third we give some combinatorial consequences of this theorem and
an illustrative example. In the first section some necessary background
material is recalled.

\section{Preliminaries}

In this section we review some useful notions and results. 
Let $\A$ be an additive category. We need the notion of \emph{approximations}.
Let $\X$ be an additive subcategory of $\A$, and let $A$ be an object of $\A$.
Then a map $X' \to A$ with $X' \in \X$ is called a \emph{right $\X$-approximation}
if the induced map $\Hom(\X,X') \to \Hom(\X,A)$ is an epimorphism.
There is the dual notion of a \emph{left $\X$-approximation}. These concepts were
introduced in~\cite{as}.

A map $f \colon A \to B$ in a category $\A$ is called \emph{right minimal},
if for every $g \colon A \to A$ such that $fg = f$, the map $g$ is
an isomorphism. There is the dual notion of a \emph{left minimal} map.
A right (left) approximation that is also a right (left) minimal map,
is called a \emph{minimal right (left) approximation}.

If $M, X$ and $Y$ are modules over a finite dimensional
algebra $\L$, then the \emph{reject of $M$ in $X$} is 
$$\rej_{M}X = \cap_{f \colon X \to M} \Ker f$$ and 
the \emph{trace of $M$ in $Y$} is $$t_M Y = \sum_{g \colon M \to Y} \Im g.$$
It is then clear that $\Hom(X / \rej_M X, M) \simeq \Hom(X,M)$ 
and $\Hom(M, t_M Y) \simeq \Hom(M,Y)$.

For a $\L$-module $T$, we let $\add T$ denote the full subcategory
of $\mod \L$ with objects all direct summands of direct sums of
copies of $T$. We let $\Fac T$ denote the full subcategory
of $\mod \L$ with objects all factors of
modules in $\add T$ and we let $\Sub T$ denote the full
subcategory with objects all submodules of modules in $\add T$.

\subsection{Tilting modules and torsion pairs}

In this section, let $\L$ be a, not necessarily hereditary, finite dimensional
algebra.
Then $T$ is called a tilting module in $\mod \L$ if
\begin{itemize}

\item[-]{$\pd_{\L} T \leq 1$},
\item[-]{$\Ext^1(T,T) = 0$} and
\item[-]{there is an exact sequence $0 \to \L \to T_0 \to T_1 \to 0$,
with $T_0,T_1$ in $\add T$.}
\end{itemize}

This is the original definition of tilting modules from~\cite{hr}, and
it was proved in~\cite{bon} that the third axiom can be replaced by the following: 
\begin{itemize}

\item[-]{the number of indecomposable direct summands of $T$ (up to isomorphism) is
the same as the number of simples.}
\end{itemize}

There is a dual concept of cotilting modules, and using the result of Bongartz, 
it follows that for hereditary algebras a module is a tilting
module if and only if it is a cotilting module.

Let $T$ be a tilting module. Then, the category $\T_T = \Fac T$
is closed under factors and extensions and is hence a torsion class.
The corresponding torsion-free class is $\F_T = \{X \mid \Hom(\T_T,X) =0 \}$.
It is well-known that $\T_T = \{Y \mid \Ext^1(T,Y) =0 \}$ and that
$\F_T = \{X \mid \Hom(T,X) =0 \} = \Sub (\tau T)$.

\subsection{A factor of the derived category and tilting objects}

Let $H$ be a hereditary algebra and
let $F = \tau^{-1}[1]$ and $\C = D^b(H)/F$, as in the introduction.

A crucial property of $\C$ is that it has a triangulated 
structure induced by the triangulated structure of $D^b(\mod H)$.
The following is a special case of
a theorem of Keller~\cite{k}.

\begin{thm}

The category $\C$ carries a canonical triangulated
structure, such that the canonical functor 
$$D^b(\mod H) \to \C$$ is a triangle functor.
\end{thm}

We use $[1]$ also to denote the suspension functor in $\C$.

In~\cite{bmrrt} it is proved that $\C$ has almost split triangles,
induced by the almost split triangles in $\D$. We denote 
by $\tau$ the corresponding translation functor in both categories.
Then there is an isomorphism
$$D\Hom_{\C}(B,\tau A) \simeq \Ext^1_{\C}(A,B)$$
for all objects $A$ and $B$ in $\C$. The notion of a tilting
object in $\C$ was defined in~\cite{bmrrt}, and
it is easy to see that an object $T$ in $\C$ 
is a tilting object if and only if the following holds:
$\Hom_{\C}(T, X[1]) = 0$ if and only if 
$X$ is in $\add T$. The object $T$ 
is called an exceptional object if $\Hom_{\C}(T, T[1]) = 0$, and it
is called an almost complete tilting object if in addition it has $n-1$ indecomposable direct 
summands, where $n$ is the number of simple $H$-modules.

The following summarises some results we need from~\cite{bmrrt}.
If $T$ is a module, we also denote its image in $\C$
by $T$.

\begin{prop} \label{many}

\begin{itemize}

\item[a)]{If $T$ is exceptional in $\mod H$, then
$T$ is also exceptional in $\C$.}
\item[b)]{If $T$ is a tilting module in $\mod H$, then $T$ is also a tilting object in $\C$.}
\item[c)]{If $\overline{T}$ is an almost complete tilting object in 
$\C$, then it has exactly two non-isomorphic complements 
in $\C$.}
\item[d)]{Let $\overline{T}$ be as above with two indecomposable complements $M$ and $M^{\ast}$ in
$\C$. Then there are triangles in $\C$:
$$M^{\ast} \to B \to M \to $$ and $$M \to B' \to M^{\ast} \to $$
where $B \to M$ and $B' \to M^{\ast}$ are minimal right $\add \overline{T}$-approximations.
\item[e)]{For a tilting object $T$
in $\C = D^b(\mod H) / F$, there is always a hereditary algebra $H'$
derived equivalent to $H$, such that the image of $T$ under this equivalence
is induced by a tilting module. }
} 
\end{itemize}

\end{prop}

\section{Cluster categories and cluster-tilted algebras}

In this section we show the main result describing the connection between a 
cluster category and the module theory of an associated cluster-tilted algebra.

For a hereditary algebra $H$ and a tilting module $T$, there are associated
torsion theories $(\T,\F)$ in $\mod H$, where $\T = \Fac T$,
and $(\X, \Y)$ in $\mod \End_H(T)^{\op}$, such that there are equivalences
of categories $\Hom_H(T,\ )\colon \T \to \Y$ and $\Ext^1_H(T,\ ) \colon \F \to \X$, see~\cite{bb}. 
In addition there
is an induced equivalence of derived categories (see~\cite{hap}) 
$$\RHom(\ ,\ ) \colon D^b(\mod H) \to D^b(\mod \End_H(T)^{\op}).$$ 
When $T$ is a tilting object in a cluster category, the functor $G= \Hom_{\C}(T,\ )$ behaves nicely on all of $\C$,
and is actually even dense.

\begin{prop}

The functor $G$ is full and dense.
\end{prop}

\begin{proof}

Let $\G = \End_{\C}(T)^{\op}$.
Then it is well known that $G$ induces an equivalence between the
additive categories $\add T$ and $\P(\G)$, the full subcategory of $\mod \G$
with objects the projective modules.
We first show that $G$ is dense. So let $Y$ be any module in $\mod \G$,
and let $0 \to P_1 \overset{g}{\rightarrow} P_0 \to Y \to 0$ be a minimal projective 
resolution.
Then there are modules $T_0$ and $T_1$ in $\add T$ and a map $\alpha \colon T_1 \to T_0$
such that $G(T_i) = P_i$ for $i=0,1$ and $G(\alpha) = g$.
There is a triangle
$$T_1 \overset{\alpha}{\rightarrow} T_0 \to A \to $$
in $\C$.
By applying $\Hom_{\C}(T,\ )$, we obtain
an exact sequence 
$$\Hom_{\C}(T,T_1) \to \Hom_{\C}(T,T_0) \to \Hom_{\C}(T,A) \to \Hom_{\C}(T,T_1[1]).$$
We now use that $T$ is tilting and Proposition~\ref{many} a) to conclude that 
$\Hom_{\C}(T,T_1[1]) =0$.
Hence $G(A) \simeq Y$, and the functor is dense.

Now, let $f \colon Y \to Z$ be a map between two indecomposable modules
in $\mod \G$. Choose minimal projective resolutions of $Y$ and $Z$:
$0 \to P_1 \overset{g}{\rightarrow} P_0 \to Y \to 0$ and
$0 \to Q_1 \overset{h}{\rightarrow} Q_0 \to Z \to 0$.
Then there are objects $T_i$ and $T'_i$ for $i=0,1$ such that
$G(T_i) = P_i$ and $G(T'_i) = Q_i$, and there are maps $\alpha$ and $\beta$
such that $G(\alpha) = g$ and $G(\beta) = h$.
Moreover,
there is a commutative diagram
$$
\xy
\xymatrix{
T_1 \ar[r]^{\alpha} \ar[d] & T_0 \ar[r] \ar[d] & A \ar[r] \ar@{.>}[d]^{\gamma} & {\ } \\
T'_1 \ar[r]^{\beta}  & T'_0 \ar[r] \ & B \ar[r] & {\ } 
}
\endxy
$$
where $\gamma$ is an induced map.
Then, using the same arguments as above, it is clear that $G(\gamma) = f$,
and the functor is dense.

\end{proof}

The functor $G$ is not faithful. 
We shall see that the maps that
are killed by $G$ are exactly the maps factoring through $\tau T$.
We have $G(\tau T) = \Hom_{\C}(T,\tau T) = \Hom_{\C}(T, T[1]) = 0$,
where the last equality uses that $\tau T \simeq T[1]$ in $\C$,
since $F= \tau^{-1}[1]$. Thus,  
there is an induced functor $\overline{G} \colon \C / \add (\tau T) \to \mod \G$.

\begin{thm} \label{main}

Let $T$ be a tilting object in $\C$, and $$G = \Hom_{\C}(T,\ ) \colon \C  
\to \mod \G.$$ Then the functor $\overline{G} \colon \C/ \add(\tau T) \to
\mod \G$ is an equivalence.
\end{thm}

By Proposition~\ref{many}, part e),
we can without loss of generality assume that $T$ is induced by a tilting module in
$\mod H$. 
We need the following observation.

\begin{lem} \label{help}

Let $M$ and $X$ be in $\mod H$, and assume $\Ext^1(M,M) = 0$.
Then $\Hom(\rej_M X, M) = 0$.
\end{lem}

\begin{proof}

Let $\overline{X} = \rej_M X$.
The canonical map 
$X \to X / \overline{X}$ induces an isomorphism $\Hom(X / \overline{X}, M) \to \Hom(X,M)$.
Consider the long exact sequence
$$0 \to \Hom(X / \overline{X}, M) \to \Hom(X,M) \to \Hom(\overline{X},M) \to 
\Ext^1(X / \overline{X},M).$$ 
To prove the claim it is clearly sufficient to show that
$\Ext^1(X / \overline{X},M) =0$.
If $\Hom(X, M)=0$, then $X / \overline{X} = 0$. 
So we can assume that there is a non-zero map
$X \to M$, and thus a monomorphism $X / \overline{X} \to M'$,
with $M'$ in $\add M$. Let 
$X / \overline{X} \to M''$ be the minimal left $\add M$-approximation,
which is then automatically a monomorphism. It induces
an epimorphism $\Ext^1(M'',M) \to \Ext^1(X / \overline{X}, M) \to 0$, and hence the claim follows
from our assumption that $\Ext^1(M,M) = 0$.
\end{proof}

We can now proceed to the proof of the theorem.

\begin{proof}

We only need to show that the functor $\overline{G}$ is faithful.
Let $f \colon A \to B$ be a map between indecomposable objects in  
$D^b(\mod H)$, and let $\overline{f} \colon A \to B$ be the
induced map in the factor category $\C$. We need to show that if $\overline{f} \neq 0$
and $G(\overline{f}) = 0$, then $f$ factors through $\add (\tau T)$.
Without loss of generality we may assume that $A$ either is a module or of the form
$P[1]$ for a projective module $P$.
We need to discuss three cases:
\begin{itemize} 

\item[I.]{Both $A$ and $B$ are in $\mod H$.}
\item[II.]{Only $A$ is in $\mod H$.}
\item[III.]{$A$ is not in $\mod H$.}
\end{itemize}

\noindent \emph{Case I.} Let $\overline{A}= A / t_T A$ and $\overline{B} = 
\rej_{\tau^2 T} B$.
We will show that there is a commutative diagram $\nabla$
$$
\xy
\xymatrix{A \ar[rr]^{f} \ar[d] & & B \\
\overline{A} \ar[rr] \ar[dr] & & \overline{B} \ar[u] \\
& C \ar[ur] & 
}
\endxy
$$
with $C$ in $\add (\tau T)$.

Our first claim is that $\Im f \subset \overline{B}$. 
Note that $\Hom_{\C}(T,A) = \amalg_i \Hom_{\D}(F^i T, A) 
= \Hom_{\D}(T,A) \ \amalg \ \Hom_{\D}(F^{-1} T,A)$,
when $T$ and $A$ are modules.
Since $G(f) = \amalg_i \Hom_{\D}(F^i T,f) = 0$ where $F = \tau^{-1}[1]$,
in particular the map $$\Hom_{\D}(\tau T [-1], A) \to 
\Hom_{\D}(\tau T [-1], B)$$ is the zero-map. 
We have a commutative diagram
$$
\xy
\xymatrix{\Hom(\tau T,A[1]) \ar[r]^0 \ar@{=}[d] & \Hom(\tau T,B[1]) \ar@{=}[d] \\
 \Ext^1(\tau T,A) \ar[r]^0 \ar@{=}[d] & \Ext^1(\tau T,B)  \ar@{=}[d]\\
D\Hom(A, \tau^2 T) \ar[r]^0  & D \Hom(B,\tau^2 T)
}
\endxy
$$
Thus, we conclude that for any map $B \to \tau^2 T$, the composition
$A \to B \to \tau^2 T$ is 0, in other words $\Im f \subset \overline{B} = \rej_{\tau^2 T} B$.
We also have that the map $\Hom(T,A) \to \Hom(T,B)$ is the zero-map,
and so it follows that $f$ factors through $\overline{A}$.

Next, we show that the induced map $\overline{A} \to \overline{B}$
factors through an object $C$ in $\add \tau T$.
First note that $\Hom(\overline{B}, \tau^2 T) =0$.
This follows from Lemma~\ref{help}, since $\Ext^1(\tau^2 T, \tau^2 T) =0$.
Since $\overline{A}$ is in $\F_T = \Sub(\tau T)$, there is an exact
sequence 
\begin{equation}\label{exactseq} 
0 \to \overline{A} \to C \to K \to 0
\end{equation}

where the map
$\overline{A} \to C$ is a minimal left $\add (\tau T)$-approximation.
We want to show that also $K$ is in $\F_T$, so we need that $\Hom(T,K) = 0$.
If we apply $\Hom(\ , \tau T)$ to~(\ref{exactseq}), we get an exact sequence
$$0 \to \Hom(K,\tau T) \to \Hom(C, \tau T) \to \Hom(\overline{A}, \tau T) 
\to \Ext^1(K,\tau T) \to \Ext^1(C,\tau T).$$
By using that the map $\overline{A} \to C$ is a minimal 
left $\add \tau T$-approximation and that $\Ext^1(\tau T, \tau T)  =0$ and hence 
$\Ext^1(C, \tau T)  =0$,
it follows that  $\Ext^1(K, \tau T) = 0$. Then we obtain 
$\Hom(T,K) = 0$, and hence $K$ is in $\F_T = \Sub(\tau T)$.
Using this and the fact that 
$\Hom(\tau^{-1}\overline{B}, \tau T) = \Hom(\overline{B}, \tau^2 T)= 0$, also
$\Hom(\tau^{-1} \overline{B},K) = 0$, and hence $\Ext^1(K,\overline{B}) = 0$.
From~(\ref{exactseq}) we now get an exact sequence
$$0 \to \Hom(K,\overline{B}) \to \Hom(C,\overline{B}) \to \Hom(\overline{A}, \overline{B}) 
\to 0.$$
Hence any map $\overline{A} \to \overline{B}$ factors through $\overline{A} \to C$, and
since $C$ is in $\add(\tau T)$, we have completed the diagram $\nabla$,
and finished case I. \\
\\
\noindent \emph{Case II.} Assume now that $A$ is in $\mod H$, while $B$ is not a module.
Since $H$ is hereditary and $\Hom(A,B) \neq 0$, we then must have
that $B = B'[1]$ for some indecomposable module $B'$. 
It is clear that $\overline{A}$ is in in $\F_T$, and that $f$ factors through
$\overline{A}$. As in case I, we can choose a minimal left
$\add \tau T$-approximation $\overline{A} \to X$, and complete
it to an exact sequence 
$$0 \to \overline{A} \to X \to K \to 0.$$     
We now apply $\Hom(\ ,B)$ to this sequence and obtain that
$$0 \to \Hom(K,B) \to \Hom(X,B) \to \Hom(\overline{A}, B) 
\to 0$$ is exact, since obviously $\Ext^1(K,B) = \Ext^1(K, B'[1]) = \Ext^2(K,B') = 0$.
This proves the claim for case II.\\
\\
\noindent \emph{Case III} Assume $A = P[1]$ for some 
indecomposable projective module $P$. As in
case II, it is clear that $B= B'[1]$ for some module $B'$.
Since $G(f) =0$, we have that in $\D$, the composition $T \to P[1] \to B'[1]$ 
vanishes for any map $T \to P[1]$, and
$\tau^{-1} T \to P \to B'$ vanishes for any map $\tau^{-1}T \to P$.
Using that
$$
\xy
\xymatrix{\Hom(T,P[1]) \ar[r]^0 \ar@{=}[d] & \Hom(T,B'[1]) \ar@{=}[d] \\
 \Ext^1(T,P) \ar[r]^0 \ar@{=}[d] & \Ext^1(T,B')  \ar@{=}[d]\\
D\Hom(P, \tau T) \ar[r]^0  & D\Hom(B',\tau T)
}
\endxy
$$
we also have that the composition $P \to B' \to \tau T$ vanishes
for any map $B' \to \tau T$.
We shall find a $C$ in $\add T$ such that there is a commutative diagram
$\nabla'$
$$
\xy
\xymatrix{P \ar[rr] \ar[d] & & B' \\
\overline{P} \ar[rr] \ar[dr] & & \overline{B'} \ar[u] \\
& C \ar[ur] & 
}
\endxy
$$
where $\overline{P} = P / t_{\tau^{-1}T}P$ and $\overline{B'} = \rej_{\tau T} B'$.
The above shows that $\Im(P \to B') \subseteq \overline{B'}$ and since
$\Hom(\tau^{-1}T,P) \to \Hom(\tau^{-1}T,B')$ is the zero-map,
the map $P \to \overline{B'}$ factors through $\overline{P}$.
We have $\Hom(\overline{B'}, \tau T) = 0$, from Lemma~\ref{help}.
We claim that $\overline{P}$ is in $\Sub T$. Since $T$ is also a cotilting module,
we have $\Sub T = \{X \mid \Ext^1(X,T) =0 \}$. 
Thus, we need to show that $\Hom(\tau^{-1} T, \overline{P})=0$.
Let $tP = t_{\tau^{-1}T}P$.
Since $tP$ is a factor of a finite direct sum of copies of the exceptional module $\tau^{-1}T$,
it is clear that $\Ext^1(\tau^{-1}T, tP) = 0$.
There is a long exact sequence
$$0 \to \Hom(\tau^{-1}T, tP) \simeq 
\Hom(\tau^{-1}T, P) 
\to \Hom(\tau^{-1}T, P / tP) \to
\Ext^1(\tau^{-1}T, tP)$$
and thus $\Ext^1(\overline{P},T) \simeq D\Hom(\tau^{-1}T,P/tP) = 0$,
which proves the claim that $\overline{P}$ is in $\Sub T$.

Consider the exact sequence
\begin{equation}\label{exact2} 
0 \to \overline{P} \to T' \to L \to 0
\end{equation}

where the map $\overline{P} \to T'$ is a minimal left $\add T$-approximation. 
Applying $\Hom(\ ,T)$, we have a long exact sequence 
$$\Hom(T',T) \to \Hom(\overline{P}, T) \to \Ext^1(L,T) \to \Ext^1(T' ,T),$$
where the last term vanishes, and the first map is an epimorphism.
Therefore $\Ext^1(L,T) = 0$ and $L$ is in $\Sub T$, so there is an embedding
$L \hookrightarrow T''$ with $T''$ in $\add T$.
We have that $$D\Ext^1(L,\overline{B}) \simeq \Hom(\tau^{-1} \overline{B}, L) \hookrightarrow
\Hom(\tau^{-1} \overline{B}, T'')\simeq \Hom(\overline{B}, \tau T'') = 0,$$
and hence $\Ext^1(L,\overline{B}) =0$.

We therefore get a short exact sequence
$$0 \to \Hom(L, \overline{B}) \to \Hom(T',\overline{B}) \to \Hom(\overline{P}, \overline{B}) \to 0$$
induced from~(\ref{exact2}).
Thus, every map $\overline{P} \to \overline{B}$ factors through $T'$, so we can choose
$C = T'$ to complete the diagram $\nabla'$.
We have now obtained that every map $A = P[1] \to B = B'[1]$ factors 
through $T'[1]$ with $T'$ in $\add T$. We now use 
that $F^{-1}(T[1]) = \tau T$
to finish the proof of case III, and the proof of the
theorem.
\end{proof}

The result has an especially nice consequence in the case of 
finite type.

\begin{cor}
Let $H$ be a hereditary algebra of finite representation type, and
let $T$ be a tilting object in $\C$, with $\G = \End_{\C}(T)^{\op}$. Then $\G$ has
the same number of non-isomorphic indecomposable modules as $H$.
\end{cor}

\begin{proof}
Let $h$ be the number of indecomposable modules for $H$, and let $n$ be the number of simples.
The number of indecomposable objects in $\C$ is $h +n$. The number of indecomposable
summands of $T$ and thus of $\tau T$ is $n$, and thus the number
of indecomposables for $\G$ is $h+n-n = h$.
\end{proof}

\section{Combinatorial properties}

Let $H$ be a hereditary algebra with a tilting module $T$ in $\mod H$. Let 
$\L = \End_H(T)^{\op}$ be the corresponding tilted algebra and $\G = \End_{\C}(T)^{\op}$
the cluster-tilted algebra. In this section we point out
a nice consequence of the equivalence of Theorem~\ref{main}, namely 
that the AR-quiver of $\G$ can be obtained
directly from the AR-quiver of $H$. 
We illustrate this by a small concrete example. This example also suggests
that there is a combinatorially nice relationship between the AR-quivers of
$\G$ and of $\L$, but we do not have any general result in this direction.

\subsection{The AR-quiver of $\G$}

Here we explain why and how the AR-structure of $\G$ can be obtained from the 
AR-structure of $H$.
The indecomposable objects of the derived category of a hereditary algebra are all
(isomorphic to) stalk complexes. The AR-structure of the 
derived category of $H$ is well-known. 
The AR-quiver consists of a countable number of copies of the AR-quiver of $H$, glued together using that the
translate of a projective is the corresponding injective (given by the Nakayama functor)~\cite{hap}.
The following is proved in~\cite{bmrrt}.

\begin{prop}
The category $\C$ has AR-triangles, induced by the AR-triangles in $D^b(\mod H)$.
\end{prop}

\begin{proof}
This is a special case of Proposition~1.3 in~\cite{bmrrt}.
\end{proof}

Combining this with the following observation, we see
that the AR-quiver of $\G$ can be obtained
from the AR-quiver of $H$.

\begin{prop} 
The almost split sequences in $\C / \add(\tau T) \simeq \mod \G$ are induced by
almost split triangles in $\C$. 
\end{prop}

\begin{proof}
By~\cite{ar}, right (left) almost split maps are sent to right (left) almost split maps
in the factor.
Let $M$ be a non-projective indecomposable in $\mod \G$, and
let $C$ be the corresponding object in $\C / \add(\tau T)$.
Let $A \overset{f}{\rightarrow} B \overset{g}{\rightarrow} C \to $
be the AR-triangle in $\C$, and let
$\widetilde{f} \colon \widetilde{A} \to \widetilde{B}$ and
$\widetilde{g} \colon \widetilde{B} \to C$
be the induced maps in $\mod \G$. Then it is
clear that $\tau_{\G} C = \ker \widetilde{g}$.
Since the composition $\widetilde{g} \widetilde{f}$ is the zero-map,
it follows that $\widetilde{g}$ can not be an epimorphism, and since
it is left almost split, it must be a monomorphism. It also
follows that it factors through $\ker \widetilde{g}$. But since it
is irreducible, we must have $\widetilde{A} \simeq \ker{g}$.
\end{proof}

We note that since the tilted algebra $\L$ is a factor algebra of $\G$, there is an
induced embedding $\mod \L \to \mod \G$. In the example below we see that the
AR-structure is preserved by this embedding. We do not know in what generality this holds.

\subsection{An example}

In this section we illustrate Theorem~\ref{main} and the combinatorial remarks above
by a small example. 
Let $Q$ be the quiver 
$$
\xy
\xymatrix{
1 \ar[r] & 2 \ar[r] & 3
}
\endxy
$$
and let $H = kQ$ be the path algebra, where $k$ is a field.
Then the AR-quiver of the triangulated category $D^b(H)$ looks like this:
$$
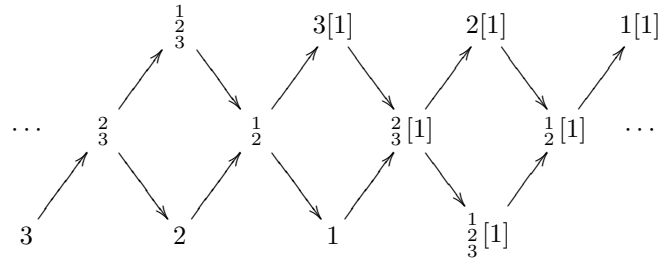
\xy
\xymatrix@R=20pt@!C=5pt{
& & {\begin{smallmatrix}1 \\2 \\ 3 \end{smallmatrix}} \ar[dr] & & 3[1] \ar[dr] & & 
2[1] \ar[dr] & & 1[1] \\
\cdots & {\begin{smallmatrix} 2 \\ 3 \end{smallmatrix}} \ar[dr] \ar[ur] & & 
{\begin{smallmatrix}1 \\2  \end{smallmatrix}} \ar[dr] \ar[ur] 
& & {\begin{smallmatrix}2 \\ 3 \end{smallmatrix}[1]} \ar[dr] \ar[ur] & & 
{\begin{smallmatrix}1 \\2 \end{smallmatrix}[1]} \ar[ur] & \cdots \\
3 \ar[ur] & & 2 \ar[ur] & & 1 \ar[ur] & & {\begin{smallmatrix}1 \\2 \\ 3 \end{smallmatrix}[1]} 
\ar[ur] & &
}
\endxy
$$

The AR-quiver of $\C$ is given by
$$
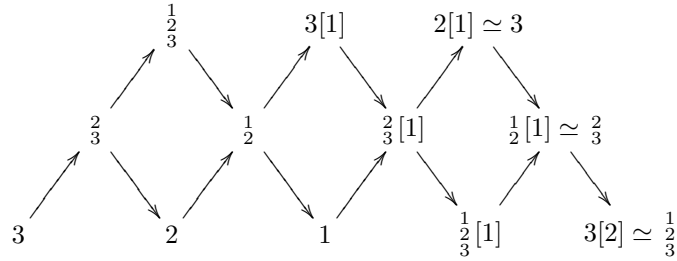
\xy
\xymatrix@R=20pt@!C=5pt{
& & {\begin{smallmatrix}1 \\2 \\ 3 \end{smallmatrix}} \ar[dr] & & 3[1] \ar[dr] & & 
2[1] \simeq 3 \ar[dr] & &  \\
& {\begin{smallmatrix} 2 \\ 3 \end{smallmatrix}} \ar[dr] \ar[ur] & & 
{\begin{smallmatrix}1 \\2  \end{smallmatrix}} \ar[dr] \ar[ur] 
& & {\begin{smallmatrix}2 \\ 3 \end{smallmatrix}[1]} 
\ar[dr] \ar[ur] & & 
{\begin{smallmatrix}1 \\2 \end{smallmatrix}[1]} \simeq  {\begin{smallmatrix} 2 \\ 3 \end{smallmatrix}}
\ar[dr] &  \\
3 \ar[ur] & & 2 \ar[ur] & & 1 \ar[ur] & & {\begin{smallmatrix}1 \\2 \\ 3 \end{smallmatrix}[1]} 
\ar[ur] & & 3[2] \simeq {\begin{smallmatrix}1 \\2 \\ 3 \end{smallmatrix}}
}
\endxy
$$
with identifications as indicated in the figure.

Let $T$ be the tilting module $T = T_1 \oplus T_2 \oplus T_3 = 
S_3 \oplus P_1 \oplus S_1$, let $\L = \End_{H}(T)^{\op}$ be the corresponding tilted
algebra and $\G =  \End_{\C}(T)^{\op}$ the cluster-tilted algebra.
We notice that $\G$
is the path algebra of the quiver

$$
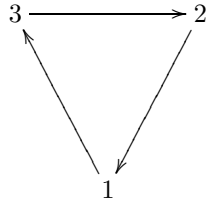
\xy
\xymatrix{3 \ar[rr] & & 2\ar[ddl] \\
& & \\
& 1 \ar[uul] &
}
\endxy
$$
with relations $\r^2=0$.
As we expect from Theorem~A, we get the AR-quiver of $\G$ by deleting the
vertices corresponding to $\tau T$.
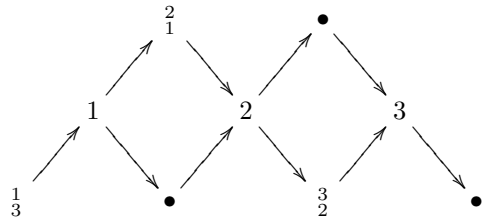
$$
\xy
\xymatrix@R=20pt@!C=5pt{
& & {\begin{smallmatrix}2 \\1 \end{smallmatrix}} \ar[dr] & & {\bullet} \ar[dr] & & \\
& 1 \ar[dr] \ar[ur] & & 2  \ar[dr] \ar[ur] 
& & 3 \ar[dr] & \\
{\begin{smallmatrix}1 \\ 3 \end{smallmatrix}} \ar[ur] & & {\bullet} \ar[ur] & & 
{\begin{smallmatrix}3 \\ 2 \end{smallmatrix}}
 \ar[ur] & & \bullet 
}
\endxy
$$

The tilted algebra $\L$ is given by the quiver $Q$ also with
$\underline{r}^2 = 0$.
We note now that there is (in this particular example) a nice embedding
of the AR-quiver of $\L$ in the AR-quiver of $\G$. The AR-quiver of
$\L$ is given by deleting the module 
$\begin{smallmatrix}1 \\ 3 \end{smallmatrix}$
in the AR-quiver above. 
This is in contrast to the embedding of $\mod \L$ into $D^b(\mod \L)$, which we consider now.
We have
an equivalence of categories $\RHom(T,\ ) \colon D^b(\mod H) \to D^b(\mod \L)$.

The image of a module $X$ in $D^b(\mod \L)$ when applying $\RHom(T,\ )$ is a complex
with homology in degree zero, given by $\Hom_{H}(T,X)$, and 
homology of degree one, given by $\Ext^1_{H}(T,X) = \Hom(T,X[1])$.
We use the numbers $1,2,3$ to indicate the simple tops of 
$\Hom(T,T_1)$, $\Hom(T,T_2)$ and $\Hom(T,T_3)$, respectively.

$$
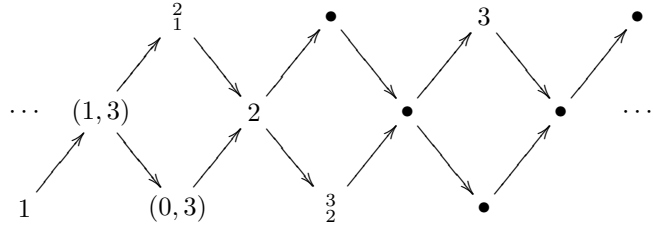
\xy
\xymatrix@R=20pt@!C=5pt{
& & {\begin{smallmatrix}2 \\1 \end{smallmatrix}} \ar[dr] & & {\bullet} \ar[dr] & & 3 \ar[dr] & & 
{\bullet} \\
\cdots & {(1,3)} \ar[dr] \ar[ur] & & 2  \ar[dr] \ar[ur] 
& & {\bullet} \ar[dr] \ar[ur] & & {\bullet} \ar[ur] &  \cdots \\
1 \ar[ur] & & (0,3) \ar[ur] & & {\begin{smallmatrix}3 \\ 2 \end{smallmatrix}} \ar[ur] & & {\bullet} 
\ar[ur] &&
}
\endxy
$$
The pair $(a,b)$ indicates that the homology in degree zero is $a$, and
the homology in degree one is $b$.
For stalk complexes in degree 0 we write $a$ instead of $(a,0)$.

Now, we do \emph{not} get the AR-quiver for $\mod \L$ directly when restricting
from $D^b(\L)$ to $\mod \L$,
but have to move some of the modules, to take the place of complexes to get
the AR-quiver of $\mod \L$.
The problem is that we normally will have some non-projective $\L$-module
$Y$ for which $\tau_{D^b(\L)} Y$ is not a $\L$-module. So when we restrict, and
$\tau Y$ gets deleted, there is then no  $\tau_{D^b(\L)} Y$-correspondent 
for $Y$, without moving.

\subsection{Some remarks}

We observe that 
the endomorphism-ring of a tilting object in $\C$ can both
have cycles in the AR-quiver and have infinite global dimension.
So in this respect, the algebras that appear as cluster-tilted algebras
are very different from the ordinary tilted algebras.

We also make some comparison with the combinatorial setup
in the classification of selfinjective algebras of finite type~\cite{rie}.
To find the possible AR-quivers for a selfinjective algebra
of finite type, it was first proved that the possible stable AR-quivers 
were given by $\Z \Delta / G$ for a Dynkin diagram $\Delta$ and $G$ a group acting 
admissibly. Then one ``filled in'' extra vertices 
for the projectives to get a module category.
Here the situation is ``opposite''. One has the translation quiver
of $\C$, and then removes certain vertices 
(given by $\tau T$) to get the translation quiver for a module category.

It is an interesting phenomenon that one can remove 
a set of vertices from a stable translation quiver, like the AR-quiver of $\D/F$,
and then obtain an AR-quiver. A necessary condition for this is 
that for any indecomposable projective $P$, we have $\tau^2 P \simeq I'$
for an indecomposable injective $I'$.

\section{Generalised APR-tilting}

Given a tilting module $T$ for a hereditary algebra $H$, 
we have shown that there is an equivalence
$\C / \add (\tau T) \to \mod \G$, where $\G = \End_{\C}(T)^{\op}$.
We now apply this to get a generalisation of the equivalence coming from 
APR-tilts.

Let $\overline{T}$ be an almost complete tilting object in $\C$, with complements
$M$ and $M^{\ast}$. By Proposition~\ref{many}, there are
triangles $M^{\ast} \to B \to M \to $ and $M \to B' \to M^{\ast} \to$,
with $B$ and $B'$ in $\add \overline{T}$.
Let $T = \overline{T} \oplus M$, and let $T' = \overline{T} \oplus M^{\ast}$. Furthermore,
let $\G = \End_{\C}(T)^{\op}$ and $\G' = \End_{\C}(T')^{\op}$. 
Then $\Hom_{\C}(T,M)$ and $\Hom_{\C}(T',M^{\ast})$ are indecomposable projective modules 
over $\G$ and $\G'$, respectively.
Thus, they have simple tops, which we denote by $S_M$ and $S_{M^{\ast}}$.
The \emph{top} of a module $Z$ is $Z/\underline{r}Z$, where $\underline{r}$ is the radical of
the ring.
We now have an exact sequence
\begin{equation}\label{long} 
\Hom_{\C}(T,B) \to \Hom_{\C}(T,M) \to \Hom_{\C}(T,\tau M^{\ast}) \to \Hom_{\C}(T,B[1])
\end{equation}

where the last term vanishes. 
\begin{lem}
The $\G$-module $\Hom_{\C}(T,\tau M^{\ast})$ is simple.
\end{lem}

\begin{proof}
By~\cite{bmrrt}, $\Hom_{\C}(M, M^{\ast}[1])$ is
is one-dimensional over the factor algebra $\End_{\C}(M)/ \Rad(M,M)$, and thus 
a simple $\End_{\C}(M)$-module. We have 
\begin{multline*}
\Hom_{\C}(T, \tau M^{\ast}) = \Hom_{\C}(\overline{T} \oplus M, \tau M^{\ast}) \simeq \\
\Hom_{\C}(\overline{T}, M^{\ast}[1]) \oplus \Hom_{\C}(M, M^{\ast}[1]) \simeq 
\Hom_{\C}(M, M^{\ast}[1]). \end{multline*}
Thus, $\Hom_{\C}(T,\tau M^{\ast}) \simeq \Hom_{\C}(M, M^{\ast}[1])$ as a $\G$-module, and
is hence simple.
\end{proof}

From this lemma and the exact sequence~(\ref{long}), it follows that
$S_M \simeq \Hom(T,\tau M^{\ast})$. Similarly we get
$S_{M^{\ast}}\simeq \Hom(T', \tau M)$.
We can apply the result from the previous section.
By Theorem~\ref{main} we get an equivalence $\C/ \add \tau T \to \mod \G$ such
that $\tau M^{\ast} \mapsto S_M$.
We also get an equivalence  $\C/ \add \tau T' \to \mod \G'$
such that $\tau M \mapsto S_{M^{\ast}}$.
If we now let $\widetilde{T} = \overline{T} \oplus M \oplus M^{\ast}$, we get equivalences
$\C / \add \widetilde{T} \to \mod \Gamma / \add S_M$ and
$\C / \add \widetilde{T} \to \mod \Gamma ' / \add S_{M^{\ast}}$

Putting this together, we get the following.

\begin{thm} \label{main2}
Let $\overline{T}$ be an almost complete tilting object in $\C$ with
complements $M$ and $ M^{\ast}$. 
Let $T = \overline{T} \oplus M$, let $T' = \overline{T} \oplus M^{\ast}$ and let
$\G = \End_{\C}(T)^{\op}$ and $\G' = \End_{\C}(T')^{\op}$.
Furthermore, let $S_M$ and $S_{M^{\ast}}$ denote the simple tops of
$\Hom_{\C}(T,M)$ and $\Hom_{\C}(T',M^{\ast})$, respectively. Then there is 
an equivalence $\mod \G / \add S_M \to \mod \G'/ \add S_{M^{\ast}}$. 
\end{thm}

Consider the following special case.
Let $M$ be an indecomposable projective $H$-module, and assume $H$ as a
left $H$-module
decomposes into $P \oplus M$, where $M$ is not in $\add P$.
Then $P$ is an almost complete tilting module.
As before, denote by $M^{\ast}$ the second complement to $P$ in $\C$.
Now, we have $\G = \End_{\C}(H)^{\op} = \End_{H}(H)^{\op} \simeq H$.
With $\G' = \End_{\C}(P \oplus M^{\ast})$,
the above theorem now says that there is an equivalence 
$\mod H / \add S_M \to \mod \G' / \add S_{M^{\ast}}$. \\
\\
{\bf Example:} Consider again the Example in
Section 3.2. With $P = S_3 \oplus P_1$, we have
the complements $M = P_2$ and $M^{\ast} = S_1$.
Now we adopt the notation in Section 3.2, but denote the cluster-tilted
algebra by $\G'$ (i.e. $\G'=\End_{\C}(P\oplus M^{\ast})^{\op}$).
Then, we have in fact an equivalence
$$\mod kQ / \add S_M \simeq \mod \G' / \add S_{M^{\ast}},$$
where $S_M$ is the simple $kQ$-module corresponding to vertex 2 in $Q$,
and $S_{M^{\ast}}$ is the simple $\G'$-module corresponding to vertex 2 in
the quiver of $\G'$. \\ 
\\
\\
Return now to the situation with $M$ indecomposable projective.
If we in addition assume that $M=S$ is simple, we are in the classical
APR-tilting setting~\cite{apr}.
It is easy to see that in this case $M^{\ast} = \tau^{-1}S$ is the 
second complement to $P$.
The module $T' = P \oplus \tau^{-1} S$ is called an APR-tilting module and
we have $\G' = \End_{\C}(T)^{\op} \simeq \End_H(T)^{\op}$.
Theorem~\ref{main2} now states that
$\mod H / \add S$ and $\mod \G' / \add S_{M^{\ast}}$ are equivalent.
A reformulation of this in terms of subcategories of indecomposable modules is that
$\ind H \setminus \{S\}$ and $\ind \G \setminus \{S_{M^{\ast}} \}$ are equivalent.
This equivalence was shown in~\cite{apr}.

Interpreting this in terms of the quiver $Q$ of a hereditary 
algebra $H =kQ$, this means that while classical APR-tilting allows
us to tilt at a vertex which is a sink, {\em generalised APR-tilting}
allows us to tilt at any vertex.
The endomorphism ring of a classical APR-tilting module is 
again a hereditary algebra given as the path algebra of the 
quiver $Q'$, obtained by reversing
the arrows pointing to the vertex where we tilt. The example
of Section 3.2 gives us an indication of what to expect
for the endomorphism rings when we tilt at arbitrary vertices,
but we do not at present have a general description in terms of $Q$.

\end{document}